\documentclass[psamsfonts,reqno]{amsart}
\usepackage{amssymb,eucal,latexsym}

\theoremstyle{plain}

\newtheorem{lemma}{Lemma}[section]
\newtheorem{theorem}[lemma]{Theorem}

\newtheorem{corollary}[lemma]{Corollary}

\newtheorem*{sclaim}{Claim}

\newtheorem*{stat}{\name}
\newcommand{\name}{testing}

\theoremstyle{definition}
\newtheorem{definition}[lemma]{Definition}

\theoremstyle{remark}
\newtheorem{remark}[lemma]{Remark}
\newtheorem{notation}[lemma]{Notation}

\newenvironment{all}[1]{\renewcommand{\name}{#1}\begin{stat}}
                        {\end{stat}}

\newcommand{\qedc}{{\qed}~{\rm Claim~{\theclaim}.}}
\newcommand{\qedsc}{{\qed}~{\rm Claim.}}

\newenvironment{scproof}
{\begin{proof}[Proof of Claim.]}
{\qedsc\renewcommand{\qed}{}\end{proof}}

\numberwithin{equation}{section}

\newcommand{\eps}{\varepsilon}
\newcommand{\ol}[1]{\overline{#1}}
\newcommand{\les}{\leqslant}
\newcommand{\pure}{\les_{\mathrm{pure}}}
\newcommand{\subd}{\les^{\oplus}}

\newcommand{\AC}[2]{G_{#1}[#2]}

\newcommand{\pup}[1]{\textup{(}{#1}\textup{)}}

\newcommand{\Mat}{\mathrm{M}}
\newcommand{\NN}{\mathbb{N}}
\newcommand{\ZZ}{\mathbb{Z}}

\newcommand{\HH}{\mathcal{H}}

\newcommand{\into}{\hookrightarrow}
\newcommand{\onto}{\twoheadrightarrow}
\newcommand{\znzo}{(\ZZ/n\ZZ)\sqcup\set{0}}

\newcommand{\jirr}{join-ir\-re\-duc\-i\-ble}

\newcommand{\res}{\mathbin{\restriction}}

\newcommand{\set}[1]{\left\{#1\right\}}
\newcommand{\setm}[2]{\set{{#1}\mid{#2}}}
\newcommand{\famm}[2]{\left({#1}\mid{#2}\right)}

\DeclareMathOperator{\J}{J} \DeclareMathOperator{\Id}{Id}
\DeclareMathOperator{\Sub}{Sub}

\newcommand{\bD}{\boldsymbol{D}}

\newcommand{\cB}{\mathcal{B}}

\newcommand{\cL}{\mathcal{L}}
\newcommand{\Rep}{\mathcal{R}_{\mathrm{ep}}}
\newcommand{\oB}{\ol{\mathcal{B}}}
\newcommand{\oL}{\ol{\mathcal{L}}}
\newcommand{\oR}{\ol{\mathcal{R}}}
\newcommand{\cO}{\mathcal{O}}
\newcommand{\cOI}{\cO_{\!\infty\,}}

\newcommand{\mbar}{\overline{m}}

\begin{document}

\author[E. Pardo]{Enrique Pardo}
\address{Departamento de Matem\'aticas\\
Universidad de C\'adiz\\
Apartado 40\\
11510 Puerto Real (C\'adiz)\\
Spain} \email{enrique.pardo@uca.es}
\urladdr{http://www2.uca.es/dept/matematicas/PPersonales/PardoEspino/EMAIN.HTML}

\author[F. Wehrung]{Friedrich Wehrung}
\address{CNRS, UMR 6139\\
D\'epartement de Math\'ematiques\\
Universit\'e de Caen\\
14032 Caen Cedex\\
France} \email{wehrung@math.unicaen.fr}
\urladdr{http://www.math.unicaen.fr/\~{}wehrung}

\title[Semilattices of groups]{Semilattices of groups and
nonstable K-theory of extended Cuntz limits}

\thanks{The research of the first author was partially supported
by the DGI and European Regional Development Fund, jointly,
through Project MTM2004-00149, by PAI III grant FQM-298 of the
Junta de Andaluc\'{\i}a, and by the Comissionat per Universitats
i
Recerca de la Generalitat de Catalunya.}


\begin{abstract}
We give an elementary characterization of those abelian
monoids~$M$ that are direct limits of countable sequences of
finite direct sums of monoids of the form either
$(\ZZ/n\ZZ)\sqcup\set{0}$ or $\ZZ\sqcup\set{0}$. This
characterization involves the Riesz refinement property together
with lattice-theoretical properties of the collection of all
subgroups of $M$ (viewed as a semigroup),
and it makes it possible to express~$M$ as a certain
submonoid of a direct product $\Lambda\times G$, where~$\Lambda$
is a distributive semilattice with zero and~$G$ is an abelian
group. When applied to the monoids $V(A)$ appearing in the
nonstable K-theory of C*-algebras, our results yield a full
description of $V(A)$ for C*-inductive limits $A$ of finite sums
of full matrix algebras over either Cuntz algebras $\cO_n$, where
$2\leq n<\infty$, or corners of~$\cOI$ by projections, thus
extending to the case including~$\cOI$ earlier work by the
authors
together with K.\,R. Goodearl.
\end{abstract}

\maketitle
\section{Introduction}\label{S:Intr}

The goal of this paper is the full elucidation of the nonstable
K-theory of a class of C*-algebras called \emph{extended Cuntz
limits}, defined as the C*-inductive limits of sequences of
finite
direct sums of full matrix algebras over the Cuntz
algebras~$\cO_n$ and over nonzero corners of~$\cOI$ by
projections. (We recall the definition of the latter for the
information of non-C*-algebraic readers: for $2\leq n<\infty$,
the
Cuntz algebra~$\cO_n$, introduced in~\cite{Cuntz}, is the unital
C*-algebra generated by elements $s_1$,\dots, $s_n$ with
relations
$s_i^*s_j=\delta_{ij}$ and $\sum_{i=1}^n s_is_i^*=1$. Further,
$\cOI$ is the unital C*-algebra defined by generators $s_i$,
$i\in
\NN$ and relations $s_i^*s_j=\delta_{ij}$.) Hence our work is a
continuation of \cite{GPW} (where the case of $\cOI$ is not
covered), so it provides an analogue, for extended Cuntz limits,
of the description of the range of the invariant for separable AF
C*-algebras (namely, ordered $K_0$) by Elliott~\cite{Elli} and
Effros, Handelman, and Shen~\cite{EHS}.

We begin by sketching the source of the problem and giving a
precise formulation. Most of the remainder of the paper is purely
semigroup-theoretic, except for the applications to C*-algebras
in
the final two sections.

In \cite{Ro}, R\o rdam gives a K-theoretic classification of even
Cuntz limits (i.e., C*-inductive limits of sequences of finite
direct sums of matrix algebras over $\cO_n$s with~$n$ even). The
invariant which R\o rdam used for his classification is
equivalent, in the unital case, to the pair $(V(A),[1_A])$ where
$V(A)$ denotes the (additive, commutative) monoid of Murray-von
Neumann equivalence classes of projections (self-adjoint
idempotents) in matrix algebras over a C*-al\-ge\-bra~$A$, and
$[1_A]$ is the class in $V(A)$ of the unit projection in $A$
(cf.~\cite[Sections~4.6, 5.1, and~5.2]{Black}). Thus, the unital
case of the classification states that if $A$ and $B$ are unital
even Cuntz limits, then $A\cong B$ if and only if $(V(A),[1_A])
\cong (V(B),[1_B])$, that is, there is a monoid isomorphism $V(A)
\rightarrow V(B)$ sending $[1_A]$ to $[1_B]$
(cf.~\cite[Theorem~7.1]{Ro}). Subsequently, Lin and Phillips
\cite{LPh} extended R\o rdam's classification result, by
including not only~$\cO_n$s with~$n$ even, but also nonzero corners over
$\cOI$ (i.e. \emph{extended} even Cuntz limits). While the authors
and Goodearl were writing~\cite{GPW}, R\o rdam
communicated to us \cite{Rotwo} that his classification can be
extended to all Cuntz limits by applying the work of Kirchberg
\cite{Kirch} and Phillips~\cite{Phil}. By the same reason, Lin
and Phillips' classification result can be enlarged to all extended
Cuntz limits.

Most of the paper is devoted to the proof of a
semigroup-theoretical result, namely Theorem~\ref{T:oL=oR}, that
provides an ``internal'' characterization of direct limits of
sequences of finite sums of monoids of the form either
$(\ZZ/n\ZZ)\sqcup\set{0}$ or $\ZZ\sqcup\nobreak\set{0}$. It turns
out that the hard core of the proof of Theorem~\ref{T:oL=oR}
consists of a lattice-theoretical statement about homomorphisms
from finite distributive lattices to subgroup lattices of abelian
groups, see Theorem~\ref{T:Approx}.

\section{Lattices and abelian groups}\label{S:LattGrp}

A \emph{lattice} is a structure $(L,\leq,\vee,\wedge)$, where
$(L,\leq)$ is a partially ordered set in which every pair
$\set{x,y}$ of elements admits a join, $x\vee y$, and a meet,
$x\wedge y$. The \emph{zero} (resp., \emph{unit}) of a lattice
$L$
are its smallest (resp., largest) element if it exists, then
denoted by $0$ (resp., $1$). We say that $L$ is \emph{complete},
if every subset $X$ of $L$ has a supremum, denoted by $\bigvee
X$.
For elements $a$ and $b$ in $L$, let $a\prec b$ hold, if $a<b$
and
there exists no $x\in L$ such that $a<x<b$.  A nonzero element
$p$
in $L$ is \emph{\jirr}, if~$p$ is not the join of two smaller
elements. In case $L$ is finite, there exists a largest element
of
$L$ smaller than $p$, denoted by $p_*$ (so $p_*\prec p$). We
denote by $\J(L)$ the set of all \jirr\ elements of $L$. For
$a\in
L$, we denote by $\J_L(a)$, or $\J(a)$ if $L$ is understood, the
set of all \jirr\ elements of $L$ below~$a$. It is well-known
that
if $L$ is finite, then $a=\bigvee\J(a)$ for all $a\in L$.

We say that $L$ is \emph{distributive} (resp., \emph{modular}),
if
$x\wedge(y\vee z)=(x\wedge y)\vee(x\wedge z)$ (resp., $x\geq z$
implies that $x\wedge(y\vee z)=(x\wedge y)\vee z$), for all
$x,y,z\in L$.

For an abelian group~$G$, we denote by $\Sub G$ the lattice of
all
subgroups of~$G$, ordered by inclusion. It is well-known that
$\Sub G$ is a modular lattice. For subgroups~$A$ and~$B$ of an
abelian group $G$, we abbreviate
\begin{align*}
A\les B,&\quad\text{if}\quad A\text{ is a subgroup of }B,\\
A\pure B,&\quad\text{if}\quad A\text{ is a pure subgroup of
}B,\\
A\subd B,&\quad\text{if}\quad A\text{ is a direct summand of }B.
\end{align*}
In particular, $A\subd B\Longrightarrow A\pure B\Longrightarrow
A\les B$.

We shall denote disjoint union by the symbol $\sqcup$.

\section{Equivalence of projections in
C*-algebras}\label{S:Basic}

We shall denote by $M\sim N$ the Murray-von~Neumann equivalence
of
self-adjoint, idempotent matrices $M$ and $N$ over a
C*-algebra~$A$, that is, $M\sim N$ if there exists a matrix $X$
such that $M=X^*X$ and $N=XX^*$ (in particular, $X=XX^*X$). We
denote by $[M]$ the $\sim$-equivalence class of a matrix $M$, and
by $V(A)$ the abelian monoid of $\sim$-equivalence classes of
self-adjoint, idempotent matrices over $A$, with
addition defined by $[M]+[N]=\Bigl[\begin{pmatrix}M&0\\
0&N\end{pmatrix}\Bigr]$. The monoid $V(A)$ encodes the so-called
\emph{nonstable K-theory of $A$}. We shall use the following
basic
lemma.

\begin{lemma}[folklore]\label{L:DecKThy}
Let $p$ be a projection \pup{i.e., a self-adjoint, idempotent
element} in a C*-algebra $A$ and let $\alpha,\beta\in V(A)$. If
$[p]=\alpha+\beta$, then there are projections $a,b\in pAp$ such
that $p=a+b$, $[a]=\alpha$, and $[b]=\beta$.
\end{lemma}

(Observe that the given conditions imply that $ab=ba=0$.)

\begin{proof}
Let $M\in\Mat_k(A)$ and $N\in\Mat_l(A)$ be self-adjoint,
idempotent matrices such that $\alpha=[M]$ and $\beta=[N]$. By
assumption, there exists $X\in\Mat_{k+l,1}(A)$ such that
\[
p=X^*X,\quad \begin{pmatrix}M&0\\ 0&N\end{pmatrix}=XX^*.
\]
Write $X=\begin{pmatrix}U\\ V\end{pmatrix}$ with
$U\in\Mat_{k,1}(A)$ and $V\in\Mat_{l,1}(A)$. Hence $p=U^*U+V^*V$,
$M=UU^*$, $N=VV^*$, while $UV^*=0$. {}From $XX^*X=X$ it follows
that $UU^*U=U$ and $VV^*V=V$. Therefore, $a=U^*U$ and $b=V^*V$
are
as required.
\end{proof}

\begin{lemma}[folklore]\label{L:MorEqRR0}
Let $A$ be a C*-algebra. Then:
\begin{enumerate}
\item For every $n\geq 1$, $V(\Mat_n(A))\cong V(A)$. \item If $A$
is separable, then $V(A)$ is countable. \item If $A$ is unital
and
has real rank zero, then given a nonzero projection $p\in A$,
$V(pAp)\cong V(A)|_{[p]}$ \pup{the order-ideal of $V(A)$
generated
by $[p]\in V(A)$}.
\end{enumerate}
\end{lemma}

\begin{proof}
(i), (ii). See \cite[p.~28]{Black}.

(iii). By \cite[Theorem~7.2]{AGOP}, $A$ is a unital exchange
ring.
By \cite[p.~111]{AGOP}, $V(A)|_{[p]}=V(ApA)$. Thus, as
$V(pAp)\cong V(ApA)$ by \cite[Proof of Lemma~7.3]{AF}, the result
holds.
\end{proof}

It is routine to check that for any unital C*-algebra $A$, the
class~$[1_A]$ is an order-unit in $V(A)$, and that the canonical
isomorphism $V(\Mat_m(A))\rightarrow V(A)$ sends
$[1_{\Mat_m(A)}]$
to $m[1_A]$. Observe that the isomorphism in
Lemma~\ref{L:MorEqRR0}(iii) sends $[p]\in V(pAp)$ to $[p]\in
V(ApA)$. Also observe that, by \cite[Theorem~1]{Zhang} (see also
\cite[Proposition~3.9]{BP}), every purely infinite simple
C*-algebra has real rank zero, whence Lemma~\ref{L:MorEqRR0}(iii)
applies to Cuntz algebras.

We shall also use the fact that $V(-)$ is a functor from
C*-algebras to abelian monoids that preserves finite direct sums
and inductive (direct) limits \cite[5.2.3--5.2.4]{Black}.

\section{Distributive subgroups with respect to a lattice
homomorphism} \label{S:DistrSubg}

\begin{definition}\label{D:Distrwrtf}
For lattices $D$ and $M$ and a lattice homomorphism
$\varphi\colon
D\to M$, we say that an element $a$ of $M$ is \emph{distributive
with respect to~$\varphi$}, if the map $D\to M$, $u\mapsto
a\wedge\varphi(u)$ is a lattice homomorphism.
\end{definition}

Observe that, as $\varphi$ is a lattice homomorphism, it suffices
to verify the condition
\begin{equation}\label{Eq:CondDistr}
a\wedge\varphi(x\vee
y)\leq(a\wedge\varphi(x))\vee(a\wedge\varphi(y)), \qquad\text{for
all }x,y\in D.
\end{equation}
In particular, if $D$ is finite, the unit of $D$ is the join of
all \jirr\ elements of~$D$, so, if $a$ is distributive with
respect to $\varphi$, we get
$a\wedge\varphi(1)=\bigvee\famm{a\wedge\varphi(p)}{p\in\J(D)}$.
Observe that $a$ is distributive with respect to $\varphi$ if{f}
$a\wedge\varphi(1)$ is distributive with respect to $\varphi$. We
will use the following characterization of distributive elements.

\begin{lemma}\label{L:CharDistr}
Let $D$ be a finite distributive lattice, let $M$ be a modular
lattice, let $\varphi\colon D\to M$ be a lattice homomorphism,
and
let $a\leq\varphi(1)$ in~$M$. Then $a$ is distributive with
respect to~$\varphi$ if{f} there exists a family
$\famm{a_p}{p\in\J(D)}$ of elements of~$M$ that satisfies the
following conditions:
\begin{enumerate}
\item $a_p\leq\varphi(p)$, for all $p\in\J(D)$.

\item $p\leq q$ implies that $a_p\leq a_q$, for all
$p,q\in\J(D)$.

\item $a_p\wedge\varphi(p_*)=\bigvee\famm{a_q}{q\in\J(p_*)}$, for
all $p\in\J(D)$.

\item $a=\bigvee\famm{a_p}{p\in\J(D)}$.
\end{enumerate}
\end{lemma}

\begin{proof}
Suppose first that $a$ is distributive with respect to $\varphi$
and $a\leq\varphi(1)$, and put $a_p=a\wedge\varphi(p)$, for all
$p\in\J(D)$. Observe that (i) and (ii) are trivially satisfied.
For (iii), we compute
\begin{align*}
a_p\wedge\varphi(p_*)&=a\wedge\varphi(p_*)
&&(\text{by the definition of }a_p)\\
&=\bigvee\famm{a\wedge\varphi(q)}{q\in\J(p_*)} &&(\text{because
}a\text{ is distributive with respect to
}\varphi)\\
&=\bigvee\famm{a_q}{q\in\J(p_*)}.
\end{align*}
For (iv), we compute
\begin{align*}
a&=a\wedge\varphi(1)&&(\text{because }a\leq\varphi(1))\\
&=\bigvee\famm{a\wedge\varphi(p)}{p\in\J(D)} &&(\text{because
}a\text{ is distributive with respect to
}\varphi)\\
&=\bigvee\famm{a_p}{p\in\J(D)}.
\end{align*}
Conversely, let $\famm{a_p}{p\in\J(D)}$ satisfy (i)--(iv) above
and set $a_u=\bigvee\famm{a_p}{p\in\J(u)}$, for all $u\in D$. So
$a_u\leq\varphi(u)$, for all $u\in D$. Furthermore, as~$D$ is
distributive, every \jirr\ element of $D$ is join-prime, thus the
map $\psi\colon D\to M$, $u\mapsto a_u$ is a join homomorphism.
It
also follows from condition~(iv) that $a=a_1$.

We claim that the equality $a_v\wedge\varphi(u)=a_u$ holds for
all
$u\leq v$ in $D$. As $D$ is finite, an easy induction proof
reduces the problem to the case where $u\prec v$. Let~$p$ be a
minimal element of $D$ with the property that $p\leq v$ and
$p\nleq u$. So $p$ is \jirr, $p\wedge u=p_*$, and $p\vee u=v$. We
compute
\begin{align*}
a_v\wedge\varphi(u)&=(a_u\vee a_p)\wedge\varphi(u)
&&(\text{because }\psi\text{ is a join homomorphism})\\
&=a_u\vee(a_p\wedge\varphi(u))
&&(\text{because }M\text{ is modular and }a_u\leq\varphi(u))\\
&=a_u\vee(a_p\wedge\varphi(p)\wedge\varphi(u))
&&(\text{because }a_p\leq\varphi(p))\\
&=a_u\vee(a_p\wedge\varphi(p\wedge u))
&&(\text{because }\varphi\text{ is a lattice homomorphism})\\
&=a_u\vee(a_p\wedge\varphi(p_*))
&&(\text{because }p\wedge u=p_*)\\
&=a_u\vee a_{p_*}
&&(\text{by condition~(iii)})\\
&=a_u &&(\text{because }p_*\leq u),
\end{align*}
which completes the proof of the claim. Taking $v=1$ in this
claim
yields that $\psi(u)=a_u=a\wedge\varphi(u)$, for all $u\in D$. In
particular, $\psi$ is a meet homomorphism.
\end{proof}

\begin{theorem}\label{T:DistrSbgrp}
Let $G$ be an abelian group, let $D$ be a finite distributive
lattice, let $\varphi\colon D\to\Sub G$, $u\mapsto G_u$ be a
lattice homomorphism. Then every finitely generated subgroup
of~$G_1$ is contained in some finitely generated subgroup
of~$G_1$
that is distributive with respect to~$\varphi$.
\end{theorem}

\begin{proof}
We argue by induction on $|\J(D)|$. Denote by $\HH(\varphi)$ the
set of all families $\vec A=\famm{A_p}{p\in\J(D)}$ of finitely
generated subgroups of~$G$ such that $A_p\les G_p$ and $p\leq q$
implies that $A_p\les A_q$, for all $p,q\in\J(D)$. We put
$A_u=\sum\famm{A_p}{p\in\J_D(u)}$, for all $u\in D$. Further, we
set $N(\vec A)=\setm{p\in\J(D)}{A_p\cap G_{p_*}=A_{p_*}}$. For
$\vec A,\vec B\in\HH(\varphi)$, let $\vec A\les\vec B$ hold, if
$A_p\les B_p$ for all $p\in\J(D)$.

\begin{sclaim}
For all $\vec A\in\HH(\varphi)$ and all $p\in\J(D)$, there exists
$\vec B\in\HH(\varphi)$ such that $\vec A\les\vec B$ and $N(\vec
A)\cup\set{p}\subseteq N(\vec B)$.
\end{sclaim}

\begin{scproof}
As $A_p\cap G_{p_*}\les G_{p_*}=\sum\famm{G_q}{q\in\J_D(p_*)}$
and
$A_p\cap G_{p_*}$ is finitely generated (because $A_p$ is), there
are finitely generated subgroups $H_q\les G_q$, for
$q\in\J_D(p_*)$, such that
\begin{equation}\label{Eq:DefvecH}
A_p\cap G_{p_*}\les\sum\famm{H_q}{q\in\J_D(p_*)}.
\end{equation}
As the interval $D'=\setm{x\in D}{x\leq p_*}$ is a sublattice of
$D$ with fewer \jirr\ elements (because $\J(D')$ is contained in
$\J(D)\setminus\set{p}$), it follows from the induction
hypothesis
that there exists a finitely generated subgroup~$C$ of~$G_{p_*}$
that is distributive with respect to $\varphi\res_{D'}$, and that
contains $\sum\famm{H_q}{q\in\J_D(p_*)}$. By using
\eqref{Eq:DefvecH} and the definition of $C$, we get
\begin{equation}\label{Eq:BascciC}
A_{p_*}\les A_p\cap G_{p_*}\les C\les G_{p_*}.
\end{equation}
We put
\begin{equation}\label{Eq:DefCu}
C_u=C\cap G_{p_*\wedge u}\ ,\qquad\text{for all }u\in D.
\end{equation}
By using \eqref{Eq:BascciC}, we obtain easily that
\begin{equation}\label{Eq:ContainCu}
A_{p_*\wedge u}\les C_u\les G_{p_*\wedge u},\qquad \text{for all
}u\in D.
\end{equation}
As $C$ is distributive with respect to $\varphi\res_{D'}$ and $D$
is distributive, the assignment $u\mapsto C_u$ defines a lattice
homomorphism from $D$ to $\Sub C$.

We prove that the family $\vec B=\famm{B_q}{q\in\J(D)}$, where
\[
B_q=A_q+C_q,\qquad\text{for all }q\in\J(D)
\]
(where $C_q$ is defined \emph{via} \eqref{Eq:DefCu}), is as
required for the claim. First observe that~$\vec B$ obviously
belongs to $\HH(\varphi)$. As both maps $u\mapsto A_u$ and
$u\mapsto C_u$ are join homomorphisms from $D$ to $\Sub G$, we
obtain
\begin{equation}\label{Eq:IneqBu}
B_u=A_u+C_u,\qquad\text{for all }u\in D.
\end{equation}
For all $u\geq p_*$ in $D$, it follows from \eqref{Eq:BascciC}
and
\eqref{Eq:DefCu} that $C_u=C\cap G_{p_*}=C$, and so $B_u=A_u+C$.
In particular, $B_p=A_p+C$ and (using \eqref{Eq:BascciC})
$B_{p_*}=A_{p_*}+C=C$, thus
\begin{align*}
B_p\cap G_{p_*}&=(A_p+C)\cap G_{p_*}\\
&=(A_p\cap G_{p_*})+C &&(\text{using \eqref{Eq:BascciC} and the
modularity of }\Sub
G)\\
&=C&&(\text{using \eqref{Eq:BascciC}}),
\end{align*}
so $B_p\cap G_{p_*}=B_{p_*}$, that is, $p\in N(\vec B)$.

Now let $q\in N(\vec A)$, we prove that $q\in N(\vec B)$. If
$q\leq p_*$, then, using \eqref{Eq:ContainCu} and
\eqref{Eq:IneqBu}, we get that $B_q=C_q$ and $B_{q_*}=C_{q_*}$,
so
$B_q\cap G_{q_*}=B_{q_*}$ (because~$\vec C$ belongs to
$\HH(\varphi\res_{D'})$), that is, $q\in N(\vec B)$. Suppose that
$q\nleq p_*$. {}From $p_*\wedge q<q$ it follows that $p_*\wedge
q\leq q_*$, and so $p_*\wedge q=p_*\wedge q_*$, thus (see
\eqref{Eq:DefCu}) $C_q=C_{q_*}$. It follows that
\begin{align*}
B_q\cap G_{q_*}&=(A_q+C_{q_*})\cap G_{q_*}
&&(\text{because }C_q=C_{q_*})\\
&=(A_q\cap G_{q_*})+C_{q_*} &&(\text{using \eqref{Eq:ContainCu}
and the modularity of }\Sub
G)\\
&=A_{q_*}+C_{q_*}
&&(\text{because }q\in N(\vec A))\\
&=B_{q_*} &&(\text{using \eqref{Eq:IneqBu}}),
\end{align*}
and so $q\in N(\vec B)$.
\end{scproof}

Now let $A$ be a finitely generated subgroup of $G_1$. As
$A\les\sum\famm{G_p}{p\in\J(D)}$ and $A$ is finitely generated,
there exists $\vec A=\famm{A_p}{p\in\J(D)}$ in $\HH(\varphi)$
such
that $A\les\sum\famm{A_p}{p\in\J(D)}$. By applying the Claim
above
at most $|\J(D)|$ times, starting with $\vec A$, we obtain $\vec
B=\famm{B_p}{p\in\J(D)}$ in $\HH(\varphi)$ such that $\vec
A\les\vec B$ and $N(\vec B)=\J(D)$. The latter condition means
that $B_p\cap G_{p_*}=B_{p_*}$ for all $p\in\J(D)$. Hence, by
applying Lemma~\ref{L:CharDistr} (with $M=\Sub G$), we obtain
that
the subgroup $B=\sum\famm{B_p}{p\in\J(D)}$ is distributive with
respect to $\varphi$. Furthermore, $B$ is finitely generated
(because all the $B_p$s are). Finally,
\[
A\les\sum\famm{A_p}{p\in\J(D)}\les\sum\famm{B_p}{p\in\J(D)}=B,
\]
so the subgroup $B$ is as required.
\end{proof}

\begin{definition}\label{D:PurityCond}
For a lattice $D$ and an abelian group $G$, a map $\varphi\colon
D\to\Sub G$ satisfies the \emph{purity condition}, if $u\leq v$
implies that $\varphi(u)\pure\varphi(v)$, for all $u\leq v$
in~$D$.
\end{definition}

We shall avoid the terminology ``pure homomorphism'', as it
conflicts with another one frequently used in lattice theory and
universal algebra.

\begin{remark}\label{Rk:TorFreeDistr}
In case $G$ is \emph{torsion-free}, $A\pure B$ implies that
$A\cap
C\pure B\cap C$, for all subgroups $A$, $B$, and $C$ of $G$ (this
condition is well-known to fail, as a rule, in the non
torsion-free case). In particular, in the context of
Theorem~\ref{T:DistrSbgrp}, if the original map $\varphi\colon
u\mapsto G_u$ satisfies the purity condition, then so does the
map
$u\mapsto B\cap G_u$.
\end{remark}

The lattice-theory oriented reader will observe that the proof of
Theorem~\ref{T:DistrSbgrp} depends only on a few
lattice-theoretical properties of $\Sub G$. In order to state the
corresponding lattice-theoretical generalization of
Theorem~\ref{T:DistrSbgrp}, we need the following classical
definitions. An element~$a$ in a complete lattice~$L$ is
\emph{compact}, if for every $X\subseteq L$, if $a\leq\bigvee X$,
then there exists a finite subset $Y$ of $X$ such that
$a\leq\bigvee Y$. We say that~$L$ is \emph{compactly n\oe
therian}, if it is complete, every element of~$L$ is a supremum
of
compact elements, and every subelement of a compact element
of~$L$
is compact. For example, for an abelian group $G$, the subgroup
lattice $\Sub G$ is a compactly n\oe therian modular lattice, in
which the compact elements are exactly the finitely generated
subgroups of $G$. Now we can state the announced generalization
of
Theorem~\ref{T:DistrSbgrp}. The proof is, \emph{mutatis
mutandis},
the same as the one of Theorem~\ref{T:DistrSbgrp}.

\begin{all}{Theorem~\ref{T:DistrSbgrp}$'$}
Let $D$ be a finite distributive lattice, let $M$ be a compactly
n\oe therian modular lattice, and let $\varphi\colon D\to M$ be a
lattice homomorphism. Then every compact element of $M$ below
$\varphi(1)$ lies below some compact element~$b\leq\varphi(1)$
of~$M$ such that the map $D\to M$, $u\mapsto b\wedge\varphi(u)$
defines a lattice homomorphism.
\end{all}

\section{Pure approximations of lattice homomorphisms satisfying
the purity condition}\label{S:ApproxPure}

\begin{definition}\label{D:ApproxHom}
For a lattice $D$ and an abelian group $G$, we say that a lattice
homomorphism $\varphi\colon D\to\Sub G$ satisfying the purity
condition is \emph{purely finitely approximated}, if for every
finitely generated subgroup $H$ of $G$, there exists a lattice
homomorphism $\psi\colon D\to\Sub G$ satisfying the purity
condition such that $\psi(u)$ is finitely generated and
$H\cap\varphi(u)\les\psi(u)\les\varphi(u)$, for all $u\in D$.
\end{definition}

For an abelian group $G$ and a positive integer $m$, we put
$G[m]=\setm{x\in G}{mx=0}$. We also put
$T(G)=\bigcup\famm{G[m]}{m\in\NN}$, the \emph{torsion subgroup}
of
$G$. The following lemma will make it possible to reduce the
proof
of Theorem~\ref{T:Approx} to the torsion case and the
torsion-free
case.

\begin{lemma}\label{L:DecPurLE}
Let $D$ be a lattice, let $G$ be an abelian group, and let
$\varphi\colon D\to\Sub G$, $u\mapsto G_u$ be a lattice
homomorphism satisfying the purity condition. Denote by
$\pi\colon
G\onto G/T(G)$ the canonical projection. Then each of the
following maps is a lattice homomorphism satisfying the purity
condition.
\begin{enumerate}
\item The map $\varphi[m]\colon D\to\Sub G[m]$, $u\mapsto
G_u[m]$.

\item The map $T(\varphi)\colon D\to\Sub T(G)$, $u\mapsto
T(G_u)$.

\item The map $\ol{\varphi}\colon D\to\Sub(G/T(G))$, $u\mapsto
\pi
G_u$.
\end{enumerate}
\end{lemma}

\begin{proof}
(i) (see the proof of Proposition~3.4 in \cite{GPW}). It is
obvious that $\varphi[m]$ is a meet homomorphism. Let $u,v\in D$
and let $z\in G_{u\vee v}[m]$. As $z\in G_{u\vee v}=G_u+G_v$,
there are $x\in G_u$ and $y\in G_v$ such that $z=x+y$. As
$0=mz=mx+my$, we get $mx=-my$, so $mx\in G_u\cap G_v=G_{u\wedge
v}$. As $G_{u\wedge v}\pure G_u$, there exists $t\in G_{u\wedge
v}$ such that $mx=mt$. As $z=(x-t)+(y+t)$ with $x-t\in G_u[m]$
and
$y+t\in G_v[m]$, we get $G_{u\vee v}[m]=G_u[m]+G_v[m]$.
Therefore,
$\varphi[m]$ is a lattice homomorphism.

Let $u\leq v$ in $D$, let $x\in G_v[m]$, and let $n\in\NN$ such
that $nx\in G_u[m]$. Let~$d$ be the greatest common divisor of
$m$
and $n$. There are integers $k$ and $l$ such that $km+ln=d$, so,
from $nx\in G_u[m]$ and $mx=0$ it follows that $dx\in G_u[m]$,
thus $dx\in G_u$. As $G_u\pure G_v$, there exists $y\in G_u$ such
that $dx=dy$. As~$d$ divides~$m$, we get $my=mx=0$, so $y\in
G_u[m]$. As~$d$ divides~$n$, we get $nx=ny$. Therefore,
$G_u[m]\pure G_v[m]$, thus completing the proof of (i).

As $T(G)$ is the directed union of all $G[m]$s, (ii) follows
immediately from (i).

(iii). It is obvious that $\ol{\varphi}$ is a join homomorphism.
Let $u,v\in D$ and let $\ol{x}\in\pi G_u\cap\pi G_v$, say
$\ol{x}=\pi(x)$ for some $x\in G$. There are $a,b\in T(G)$ such
that $x-a\in G_u$ and $x-b\in G_v$. Pick $m\in\NN$ such that
$ma=mb=0$. We obtain that $mx\in G_u\cap G_v=G_{u\wedge v}$,
hence, as $\varphi$ satisfies the purity condition, $mx=my$ for
some $y\in G_{u\wedge v}$, thus $x-y\in T(G)$, and so
$\ol{x}=\pi(y)\in\pi G_{u\wedge v}$. Therefore, $\ol{\varphi}$ is
a lattice homomorphism.

Let $u\leq v$ in $D$, let $\ol{x}\in\pi G_v$ and $m\in\NN$ such
that $m\ol{x}\in\pi G_u$. Writing $\ol{x}=\pi(x)$ for $x\in G_v$,
we obtain that $mx\in G_u+T(G)$, and so $nmx\in G_u$ for some
$n\in\NN$, hence, as $\varphi$ satisfies the purity condition,
$nmx=nmy$ for some $y\in G_u$, so $x-y\in T(G)$, and so
$\ol{x}=\pi(y)\in\pi G_u$. Therefore, $\ol{\varphi}$ satisfies
the
purity condition.
\end{proof}

The statements of Lemma~\ref{L:ApproxTor} and
Theorem~\ref{T:Approx} relate the concepts introduced in
Definitions~\ref{D:PurityCond} (purity condition)
and~\ref{D:ApproxHom} (purely finitely approximated). The
following result deals with the torsion case, and it is implicit
in \cite{GPW}.

\begin{lemma}\label{L:ApproxTor}
Let $D$ be a finite distributive lattice and let $G$ be an
abelian
torsion group. Then every lattice homomorphism $\varphi\colon
D\to\Sub G$ satisfying the purity condition is purely finitely
approximated.
\end{lemma}

\begin{proof}
Write $\varphi(u)=G_u$, for all $u\in D$, and let $H$ be a
finitely generated subgroup of~$G$. Pick a positive integer $m$
such that all elements of $H$ have order dividing~$m$, and put
$H_u=H\cap G_u$, for all $u\in D$. As $H\les G_1[m]$ and by
Lemma~\ref{L:DecPurLE}, this reduces the problem to the case
where
$mG=\set{0}$.

Now we argue as in the proof of \cite[Theorem~6.1]{GPW}. For all
$p\in\J(D)$, since $G_{p_*}\pure G_p$ and $mG_p=\set{0}$, it
follows from Kulikov's Theorem (see \cite[Theorem~27.5]{FuchsI})
that $G_p=G_{p_*}\oplus K_p$ for some subgroup $K_p$ of $G_p$.
Hence, \cite[Lemma~5.2]{GPW} yields that
\begin{equation}\label{Eq:GaDirSum}
G_u=G_0\oplus\bigoplus\famm{K_p}{p\in\J(u)},\quad\text{for all
}u\in D.
\end{equation}
As $H_u\les G_u$ and $H_u$ is finitely generated, for all $u\in
D$, there are finitely generated subgroups $G'_0\les G_0$ and
$K'_p\les K_p$, for $p\in\J(D)$, such that, putting
\[
G'_u=G'_0\oplus\bigoplus\famm{K'_p}{p\in\J(u)},
\]
we get $H_u\les G'_u$, for all $u\in D$. The map $u\mapsto G'_u$
is the desired approximation.
\end{proof}

Now we remove the torsion assumption from
Lemma~\ref{L:ApproxTor}.

\begin{theorem}[Pure approximation theorem]\label{T:Approx}
Let $D$ be a finite distributive lattice and let $G$ be an
abelian
group. Then every lattice homomorphism from $D$ to $\Sub G$
satisfying the purity condition is purely finitely approximated.
\end{theorem}

\begin{proof}
Let $\varphi\colon D\to\Sub G$, $u\mapsto G_u$ be a lattice
homomorphism satisfying the purity condition. Denote by
$\pi\colon
G\onto G/T(G)$ be the canonical projection and let
$\ol{\varphi}\colon D\to\Sub(G/T(G))$, $u\mapsto\pi G_u$.

Now let $A$ be a finitely generated subgroup of $G$. Without loss
of generality we may take $A\les G_1$. By applying
Theorem~\ref{T:DistrSbgrp} to the group $G/T(G)$, the
homomorphism
$\ol{\varphi}$, and the subgroup $\pi A$, we obtain a finitely
generated subgroup $\ol{H}$ of $\pi G_1$ containing $\pi A$ such
that the map $\ol{\psi}\colon D\to\Sub(G/T(G))$,
$u\mapsto\ol{H}_u=\ol{H}\cap\pi G_u$ is a lattice homomorphism.
As
$G/T(G)$ is torsion-free and $\ol{\varphi}$ satisfies the purity
condition (see Lemma~\ref{L:DecPurLE}(iii)), the map~$\ol{\psi}$
is a lattice homomorphism satisfying the purity condition (see
Remark~\ref{Rk:TorFreeDistr}). As~$\ol{H}$ is a finitely
generated
subgroup of the torsion-free abelian group $G/T(G)$, it is free
abelian (of finite rank). Denote by $\eps\colon\ol{H}\into
G/T(G)$
the inclusion map.

\begin{sclaim}
There exists a group embedding $\alpha\colon\ol{H}\into G$ such
that $\pi\circ\alpha=\eps$ and $\alpha\ol{H}_u\les G_u$ for all
$u\in D$.
\end{sclaim}

\begin{scproof}
For each $p\in\J(D)$, as $\ol{H}_{p_*}\pure\ol{H}_p$ and $\ol{H}$
is finitely generated, there exists $\ol{K}_p\les\ol{H}_p$ such
that $\ol{H}_p=\ol{H}_{p_*}\oplus\ol{K}_p$. Hence it follows from
\cite[Lemma~5.2]{GPW} that
\begin{equation}\label{Eq:DecolHu}
\ol{H}_u=\ol{H}_0\oplus\bigoplus\famm{\ol{K}_p}{p\in\J(u)},
\quad\text{for all }u\in D.
\end{equation}
For all $p\in\J(D)$, as $\ol{K}_p\les\ol{H}$ and $\ol{H}$ is free
abelian (of finite rank), $\ol{K}_p$ is free abelian (of finite
rank), thus projective. Hence, as $\ol{K}_p\les\ol{H}_p\les\pi
G_p$ and denoting by $\pi_p\colon G_p\onto\pi G_p$ the
restriction
of $\pi$ and by $\eps_p\colon\ol{K}_p\into\pi G_p$ the
restriction
of $\eps$, we obtain a group homomorphism
$\alpha_p\colon\ol{K}_p\to G_p$ such that
$\pi_p\circ\alpha_p=\eps_p$. Similarly, denoting by $\pi_0\colon
G_0\onto\pi G_0$ the restriction of $\pi$ and by
$\eps_0\colon\ol{H}_0\into\pi G_0$ the restriction of $\eps$, we
obtain a group homomorphism $\alpha_0\colon\ol{H}_0\to G_0$ such
that $\pi_0\circ\alpha_0=\eps_0$. Applying \eqref{Eq:DecolHu} to
$u=1$, we get
$\ol{H}=\ol{H}_0\oplus\bigoplus\famm{\ol{K}_p}{p\in\J(D)}$, so we
can define a group homomorphism $\alpha\colon\ol{H}\to G$ by the
rule
\begin{equation}\label{Eq:Defalpha}
\alpha\biggl(x_0+\sum_{p\in\J(D)}x_p\biggr)
=\alpha_0(x_0)+\sum_{p\in\J(D)}\alpha_p(x_p),
\end{equation}
for all $x_0\in\ol{H}_0$ and all
$\famm{x_p}{p\in\J(D)}\in\prod\famm{\ol{K}_p}{p\in\J(D)}$. {}From
$\pi_p\circ\alpha_p=\eps_p$ for all $p\in\J(D)\cup\set{0}$ it
follows that $\pi\circ\alpha=\eps$. As $\eps$ is an embedding, so
is $\alpha$. Finally, let $u\in D$ and let $x\in\ol{H}_u$. It
follows from \eqref{Eq:DecolHu} that $x$ can be decomposed as
\[
x=x_0+\sum_{p\in\J(u)}x_p,\quad\text{where }x_0\in\ol{H}_0\text{
and } x_p\in\ol{K}_p\text{ for all }p\in\J(u).
\]
As $\alpha_p(x_p)\in G_p$ for all $p\in\J(u)\cup\set{0}$, it
follows from \eqref{Eq:Defalpha} that $\alpha(x)\in G_u$.
\end{scproof}

Put $H=\alpha\ol{H}$. As $\eps$ is an embedding and
$\pi\circ\alpha=\eps$, we obtain
\begin{equation}\label{Eq:HcapTG0}
H\cap T(G)=\set{0}.
\end{equation}
As $\alpha$ is an embedding, the map $\psi\colon D\to\Sub H$,
$u\mapsto H_u=\alpha\ol{H}_u$ is a lattice homomorphism
satisfying
the purity condition (because $\ol{\psi}\colon u\mapsto\ol{H}_u$
has these properties). For $u\in D$, we observe that
\[
\pi(A\cap G_u)\les\ol{H}\cap\pi G_u=\ol{H}_u=\pi H_u,
\]
thus for all $x\in A\cap G_u$, there exists $y\in H_u$ such that
$x-y\in T(G)$. As $y\in G_u$ (because $y\in
H_u=\alpha\ol{H}_u\les
G_u$) and $x\in G_u$, we obtain that $x-y$ belongs to $T(G)\cap
G_u=T(G_u)$, and so $x\in T(G_u)+H_u$. Hence, using
\eqref{Eq:HcapTG0}, we have proved that $A\cap G_u\les
T(G_u)\oplus H_u$. As $A\cap G_u$ is finitely generated, there
exists a finitely generated subgroup $B_u$ of $T(G_u)$ such that
$A\cap G_u\les B_u\oplus H_u$.

It follows from Lemma~\ref{L:DecPurLE}(ii) that the map
$T(\varphi)\colon D\to\Sub T(G)$, $u\mapsto\nobreak T(G_u)$ is a
lattice homomorphism satisfying the purity condition. Hence,
applying Lem\-ma~\ref{L:ApproxTor} to this morphism and the sum
of
all $B_u$s, we obtain a lattice homomorphism $D\to\Sub T(G)$,
$u\mapsto G'_u$, satisfying the purity condition and with $G'_1$
finitely generated, such that $B_u\les G'_u\les T(G_u)$ for all
$u\in D$. So $A\cap G_u\les G'_u\oplus H_u$, for all $u\in D$. It
follows from \eqref{Eq:HcapTG0} that the map $D\to\Sub G$,
$u\mapsto G'_u\oplus H_u$ is a lattice homomorphism satisfying
the
purity condition, so it is as desired.
\end{proof}

As an immediate corollary, we get a lattice-theoretical
characterization of purity for embeddings of abelian groups,
similar to the one mentioned in \cite{Head}.

\begin{corollary}\label{C:Pur1Ord}
A subgroup $A$ of an abelian group $B$ is a pure subgroup if{f}
for any finitely generated $H\les B$, there are finitely
generated
$A'\les A$ and $B'\les B$ such that $A\cap H\les A'$, $H\les B'$,
and $A'\subd B'$.
\end{corollary}

\begin{proof}
That the given condition implies purity is an easy exercise (take
$H$ monogenic). Conversely, suppose that $A\pure B$ and let
$H\les
B$ be a finitely generated subgroup. Denote by $D=\set{0,1}$ the
two-element chain and by $\varphi\colon D\to\Sub B$ the
homomorphism sending $0$ to $A$ and $1$ to $B$. As $\varphi$
satisfies the purity condition and by Theorem~\ref{T:Approx},
there exists a homomorphism $\psi\colon D\to\Sub B$ with finitely
generated values such that
$H\cap\varphi(u)\les\psi(u)\les\varphi(u)$, for all $u\in D$. Put
$A'=\psi(0)$ and $B'=\psi(1)$.
\end{proof}

\section{Regular refinement monoids; the classes $\oB$, $\oL$,
and $\oR$} \label{S:DirLimExt}

We shall use the notation and terminology of \cite{GPW}
concerning
(abelian) monoids and semilattices of groups. In particular,
every
abelian monoid $M$ is endowed with a partial preordering $\leq$
defined by $x\leq y$ if{f} there exists $z$ such that $x+z=y$. We
say that $M$ is \emph{conical}, if $x+y=0$ implies that $x=y=0$,
for all $x,y\in M$. We say that $M$ is \emph{regular}, if $2x\leq
x$, for all $x\in M$, and we say that $M$ is a \emph{semilattice
of groups}, if $M$ is a disjoint union of groups (i.e.,
subsemigroups each of which happens to be a group). We say that
$M$ is a \emph{refinement monoid}, if for all $a_0,a_1,b_0,b_1\in
M$, if $a_0+a_1=b_0+b_1$, then there are $c_{i,j}\in M$, for
$i,j<2$, such that $a_i=c_{i,0}+c_{i,1}$ and
$b_i=c_{0,i}+c_{1,i}$, for all $i<2$. We put
$\Lambda(M)=\setm{x\in M}{2x=x}$. A \emph{semilattice} is an
abelian, idempotent monoid. It is \emph{distributive}, if it is a
refinement monoid.

We recall the following characterization of regular abelian
monoids, see \cite[Corollary~IV.2.2]{How}, also
\cite[Lemma~2.1]{GPW}.

\begin{lemma}\label{L:SemGrps}
An abelian monoid $M$ is regular if{f} it is a semilattice of
groups.
\end{lemma}

A regular abelian monoid $M$ is the disjoint union of all
subgroups
\[
G_M[e]=\setm{x\in M}{e\leq x\leq e},\qquad\text{for
}e\in\Lambda(M).
\]
For $a\leq b$ in $\Lambda(M)$, the assignment $j_{a,b}\colon
x\mapsto x+b$ defines the \emph{natural map} from $G_M[a]$ to
$G_M[b]$. It is a group homomorphism. If $j_{a,b}$ is an
embedding
for all $a\leq b$ in $M$, then we shall say that $M$ satisfies
the
\emph{embedding condition}, denoted by (emb). If the range of
$j_{a,b}$ is a pure subgroup of $G_M[b]$ for all $a\leq b$ in
$\Lambda(M)$, then we shall say that $M$ satisfies the
\emph{purity condition}, denoted by (pur).

The following definition introduces classes $\oB$, $\oL$, and
$\oR$ of abelian monoids, which properly contain, respectively,
the classes $\cB$, $\cL$, and $\Rep$ considered in \cite{GPW}.

\begin{notation}\label{Not:LR}
Denote by $\oB$ the class of all finite direct sums of abelian
monoids of the form $(\ZZ/n\ZZ)\sqcup\set{0}$, where $n\in\NN$,
or
$\ZZ\sqcup\set{0}$ \pup{where $\sqcup$ denotes disjoint union}.
Further, denote by $\oL$ the class of all direct limits of
monoids
from $\oB$, and by $\oR$ the class of all regular refinement
monoids satisfying the conditions (emb) and (pur).
\end{notation}

As $\oB$ is a class of finitely generated abelian monoids, closed
under finite direct sums, the following lemma is an easy
consequence of Proposition~3.1 and Section~4 in~\cite{GPW}.

\begin{lemma}\label{L:oLclosed}
The class $\oL$ is closed under direct limits, finite direct
sums,
and retracts, and contains as a member $A\sqcup\set{0}$, for any
finitely generated abelian group~$A$. Furthermore, $\oL$ is
contained in $\oR$.
\end{lemma}

\begin{lemma}\label{L:FinGenoR}
Any finitely generated monoid in $\oR$ belongs to $\oL$.
\end{lemma}

\begin{proof}
Similar to the proof of \cite[Proposition~5.3]{GPW}, the key
point
being this time that every pure subgroup of a finitely generated
abelian group is a direct summand.
\end{proof}

\begin{lemma}\label{L:DirUnFG}
Each monoid $M$ in $\oR$ is the directed union of those finitely
generated submonoids of $M$ that belong to $\oR$.
\end{lemma}

\begin{proof}
We argue as in the proof of \cite[Theorem~6.1]{GPW}. We must
prove
that any finite subset $X$ of $M$ is contained in some finitely
generated submonoid of $M$ lying in $\oR$. For convenience, we
assume that $0\in X$. Furthermore, by \cite[Theorem~3.3]{GPW}, we
may assume that
\[
M=\bigsqcup_{e\in\Lambda}(\set{e}\times
G_e)\subseteq\Lambda\times
G,
\]
for some distributive semilattice $\Lambda$ and some abelian
group
$G$ with subgroups $G_e$ satisfying the following conditions:
\begin{enumerate}
\item $G=\bigcup_{e\in\Lambda}G_e$.

\item $a\leq b$ implies that $G_a\pure G_b$, for all
$a,b\in\Lambda$.

\item $G_0=\set{0}$.

\item $G_a+G_b=G_{a+b}$, for all $a,b\in\Lambda$.

\item $G_a\cap G_b=\bigcup_{e\leq a,b\text{ in }\Lambda}G_e$, for
all $a,b\in\Lambda$.
\end{enumerate}
Now we denote by $\Id\Lambda$ the lattice of all \emph{ideals} of
$\Lambda$, that is, those nonempty subsets~$A$ of~$\Lambda$ such
that $x+y\in A$ if{f} $x,y\in A$, for all $x,y\in\Lambda$.
As~$\Lambda$ is a distributive semilattice, $\Id\Lambda$ is a
distributive lattice.

Now we set $G_A=\bigcup_{e\in A}G_e$, for all $A\in\Id\Lambda$.
Observe that the union defining~$G_A$ is directed, and that
$G_{[0,e]}=G_e$, for all $e\in\Lambda$. Hence the map
$\Id\Lambda\to\Sub G$, $A\mapsto G_A$ is a lattice homomorphism
satisfying the purity condition, and sending~$\set{0}$
to~$\set{0}$.

Write any $x\in X$ in the form $x=(e_x,g_x)\in M$. Denote by
$\bD$
the sublattice of $\Id\Lambda$ generated by $\setm{[0,e_x]}{x\in
X}$ and by $K$ the (finitely generated) subgroup of~$G$ generated
by $\setm{g_x}{x\in X}$. As $\Id\Lambda$ is distributive,~$\bD$
is
finite. Moreover, the ideal $\set{0}$ belongs to $\bD$ since
$0\in
X$. By Theorem~\ref{T:Approx}, there exists a lattice
homomorphism
$\bD\to\Sub G$, $A\mapsto G'_A$ satisfying the purity condition
such that $G'_A$ is finitely generated and $G_A\cap K\les
G'_A\les
G_A$, for all $A\in\bD$. For each $P\in\J(\bD)$, as
$G'_{P_*}\pure
G'_P$ and $G'_P$ is finitely generated, there exists $H_P\les
G'_P$ such that $G'_P=G'_{P_*}\oplus H_P$. As
$G'_{\set{0}}=\set{0}$, we obtain, using \cite[Lemma~5.2]{GPW},
that
\[
G'_A=\bigoplus\famm{H_P}{P\in\J_{\bD}(A)},\qquad \text{for all
}A\in\bD.
\]
Using the observations that $X$ is finite and that for each $x\in
X$, the element $[0,e_x]$ is the join of all \jirr\ elements of
$\bD$ below it, we obtain, as in the proof of
\cite[Theorem~6.1]{GPW}, elements $u_P\in P$, for $P\in\J(\bD)$,
such that
\[
e_x=\bigvee\famm{u_P}{P\in\J_{\bD}([0,e_x])},\qquad \text{for all
}x\in X.
\]
Since each $G'_P$ is a finitely generated subgroup of the
directed
union $G_P=\bigcup_{e\in P}G_e$, we may enlarge the elements
$u_P$
to ensure that $G'_P\les G_{u_P}$, for all $P\in\J(\bD)$.
Finally,
denoting by $P^\dagger$ the largest element of $\bD$ such that
$P\not\subseteq P^\dagger$ (its existence is ensured by $P$ being
a \jirr\ element of the finite distributive lattice~$\bD$), we
may
enlarge $u_P$ further to ensure that $u_P\in P\setminus
P^\dagger$. We define a map $\varphi\colon\bD\to\Lambda$ by the
rule
\[
\varphi(A)=\bigvee\famm{u_P}{P\in\J_{\bD}(A)},\qquad \text{for
all
}A\in\bD.
\]
Now we conclude the proof as in the final section of the proof of
\cite[Theorem~6.1]{GPW}. The construction of $\varphi$ ensures
that it is a semilattice embedding from~$\bD$ into~$\Lambda$,
that
$\varphi(A)\in A$ for all $A\in\bD$, and that
$\varphi([0,e_x])=e_x$ for all $x\in X$. Further, we set
\[
N=\bigsqcup_{A\in\bD} (\set{\varphi(A)}\times
G'_A)\subseteq\Lambda\times G.
\]
As $A\mapsto G'_A$ is a zero-preserving lattice homomorphism
satisfying the purity condition, it follows from
\cite[Theorem~3.3]{GPW} that $N$ belongs to $\oR$. As all groups
$G'_A$ are finitely generated, $N$ is finitely generated. As
\[
G'_A=\sum_{P\in\J_{\bD}(A)}G'_P\les\sum_{P\in\J_{\bD}(A)}G_{u_P}
\les\sum_{P\in\J_{\bD}(A)}G_{\varphi(P)}=G_{\varphi(A)},
\]
for all $A\in\bD$, we see that $N$ is contained in $M$. Finally,
for each $x\in X$, as~$g_x$ belongs to $G_{e_x}=G_{[0,e_x]}$ and
$g_x\in K$, we get $g_x\in G_{[0,e_x]}\cap K\subseteq
G'_{[0,e_x]}$, so $x=(e_x,g_x)\in N$. Therefore, $X$ is contained
in $N$.
\end{proof}

By using Lemmas~\ref{L:oLclosed}, \ref{L:FinGenoR}, and
\ref{L:DirUnFG}, we obtain our main monoid-theoretical result.

\begin{theorem}\label{T:oL=oR}
The direct limits of finite direct sums of abelian monoids of the
form\linebreak $(\ZZ/n\ZZ)\sqcup\set{0}$, where $n\in\NN$, or
$\ZZ\sqcup\set{0}$, are exactly the regular conical refinement
monoids satisfying \pup{emb} and \pup{pur}. That is, $\oL=\oR$.
\end{theorem}

An obvious adaptation of \cite[Corollary 6.6]{GPW} give us the
following result.

\begin{corollary}\label{C:Limorderunit}
Let $(M,u)$ be an abelian monoid with order-unit. Then $(M,u)$ is
a direct limit of finite direct sums of pairs of the form
$(\znzo,\, \mbar)$ and $\bigl( \ZZ\sqcup\set{0},\, m \bigr)$ if
and only if $M$ is a regular conical refinement monoid satisfying
\pup{emb} and \pup{pur}.
\end{corollary}

\section{Lifting monoid maps by C*-algebra
maps}\label{S:C*applics}

\begin{definition}\label{D:ExtCuntz}
A \emph{Cuntz algebra} is an algebra of the form $\cO_n$, where
$2\leq n\leq\infty$. A \emph{special Cuntz limit} is a
C*-inductive limit of a sequence of finite direct sums of Cuntz
algebras. An \emph{extended Cuntz limit} is a C*-inductive limit
of a sequence of finite direct sums of full matrix algebras over
Cuntz algebras $\cO_n$ for $2\leq n<\infty$ and nonzero corners
of
$\cOI$ by projections.
\end{definition}

The basic K-theoretic information concerning the Cuntz algebras
$\cO_n$, where $2\leq n<\infty$, is usually summarized in the
statements $K_0(\cO_n)\cong\ZZ/(n-1)\ZZ$ and $K_1(\cO_n)=0$
\cite[Theorems~3.7--3.8]{Cuntz2}. Also, $K_0(\cOI)\cong\ZZ$ and
$K_1(\cOI)=0$ \cite[Corollary 3.11]{Cuntz2}. However, Cuntz also
showed that the Murray-von Neumann equivalence classes of nonzero
projections in a purely infinite simple C*-algebra $A$ form a
subgroup of $V(A)$ that maps isomorphically onto $K_0(A)$ under
the natural map $V(A) \rightarrow K_0(A)$ \cite[p.~188]{Cuntz2}.
Moreover \cite[Proposition~2.1, Corollary~2.2]{AGP},
$V(A)\cong\set{0}\sqcup K_0(A)$. It follows that
$V(\cO_n)\setminus\set{0}$ is a group isomorphic to $K_0(\cO_n)$,
that is, $V(\cO_n)\cong (\ZZ/(n-1)\ZZ)\sqcup\set{0}$. It is
routine to check that the corresponding isomorphism sends
$[1_{\cO_n}]$ to the coset $\overline{1}$ in $\ZZ/(n-1)\ZZ$, and
thus we get
\begin{equation}\label{Eq:VMmOn}
\bigl( V(\Mat_m(\cO_n)),\, [1_{\Mat_m(\cO_n)}] \bigr) \cong
\bigl((\ZZ/(n-1)\ZZ)\sqcup\set{0},\, \mbar\bigr)
\end{equation}
for all $m\geq 1$ and $n\geq 2$.

\begin{remark}\label{Rk:Classification}
Because of \eqref{Eq:VMmOn}, for every $n\geq 2$, $k\geq 1$, and
for every nonzero projection $p\in \Mat_k(\cO_n)$, there exists
$l\in\set{1,2,\dots,n-1}$ such that $p\Mat_k(\cO_n)p\cong
\Mat_l(\cO_n)$. To see this, observe that as $p\neq0$,
there exists
$l\in\set{1,2,\dots,n-1}$ such that
$[p]=\overline{l}\in\ZZ/(n-1)\ZZ$. Denote by $\mathbb{K}$ the
C*-algebra of compact operators over an infinite-dimensional,
separable Hilbert space $\mathcal{H}$, set
$I_t=\sum\limits_{i=1}^{t}e_{i,i}\in \mathbb{K}$ for any
$t\geq1$,
and set $E_l=1\otimes I_l\in \cO_n\otimes \mathbb{K}$. {}From
\[
\Mat_k(\cO_n)\otimes\mathbb{K}
\cong\cO_n\otimes\Mat_k(\mathbb{C})\otimes\mathbb{K} \cong
\cO_n\otimes \Mat_k(\mathbb{K})
\]
and $\Mat_k(\mathbb{K})\cong\mathbb{K}$ (see
\cite[Proposition~1.10.2]{W-O}), we obtain a natural isomorphism
$\Mat_k(\cO_n)\otimes\mathbb{K}\cong \cO_n\otimes \mathbb{K}$.
Let
the projection $q\in\cO_n\otimes \mathbb{K}$ correspond to
$p\otimes I_1\in\Mat_k(\cO_n)\otimes \mathbb{K}$ under this
isomorphism (whence $[q]=\overline{l}\in \ZZ/(n-1)\ZZ$). It
follows that $q\sim E_l$ in $\cO_n\otimes \mathbb{K}$, whence
\[
p\Mat_k(\cO_n)p\cong (p\otimes I_1)(\Mat_k(\cO_n)\otimes
\mathbb{K})(p\otimes I_1)\cong q(\cO_n\otimes \mathbb{K})q\cong
E_l(\cO_n\otimes \mathbb{K})E_l\cong \Mat_l(\cO_n),
\]
as desired.
\end{remark}

Analogously \cite[Section~3]{Cuntz2}, we get $V(\cOI)\cong
\ZZ\sqcup\set{0}$, \emph{via} an isomorphism sending $[1_{\cOI}]$
to $1$ in $\ZZ$, and
\begin{equation}\label{Eq:VMmOn2}
\bigl( V(\Mat_m(\cOI)),\, [1_{\Mat_m(\cOI)}] \bigr) \cong \bigl(
\ZZ\sqcup\set{0},\, m \bigr),\qquad\text{for all }m\in\NN.
\end{equation}

We also need to consider the case of the pairs $(\ZZ, -m)$, with
$m\in \ZZ^+$, which cannot be represented by any pair
$(K_0(\Mat_n(\cOI)), [1_{\Mat_n(\cOI)}])$. We can solve this
problem by using corners by projections of $\cOI$. Throughout
Sections~\ref{S:C*applics} and~\ref{S:OrdUnits}, we shall use the
projections of $\cOI$ defined as
\[
p_n=1-\sum_{i=1}^ns_is_i^*,\text{ for all }n\geq0\quad(\text{so
}p_0=1).
\]
Observe that $[p_n]\in K_0(\cOI)$ corresponds to $-(n-1)\in \ZZ$.
Hence
\begin{equation}\label{Eq:VMmOn3}
\bigl( V(p_n\cOI p_n),\, [p_n] \bigr) \cong \bigl(
\ZZ\sqcup\set{0},\, -(n-1) \bigr)
\end{equation}

\begin{remark}\label{Rk:Classification2}
Because of \eqref{Eq:VMmOn2} and \eqref{Eq:VMmOn3}, we get the
following facts:

\noindent (i) For every $k\geq 1$ there exists a projection $p\in
\cOI$ such that $p\cOI p\cong \Mat_k(\cOI)$. To see this, notice
that for every $k\geq 1$ there exists a projection $p\in \cOI$
such that $[p]=k\in \ZZ$. So, as in
Remark~\ref{Rk:Classification}, $p\otimes I_1\sim 1\otimes I_k$
in
$\cOI\otimes \mathbb{K}$, and hence
\[
p\cOI p\cong (p\otimes I_1)(\cOI\otimes \mathbb{K})(p\otimes
I_1)\cong (1\otimes I_k)(\cOI\otimes \mathbb{K})(1\otimes
I_k)\cong \Mat_k(\cOI).
\]

\noindent (ii) For every projection $p\in \cOI$ and every
$k\in\ZZ^+$, one of the following cases occurs (by arguments
similar to those in part (i)):
\begin{itemize}
\item[(a)] If $\bigl( V(p\cOI
p),\,[p]\bigr)\cong\bigl(\ZZ\sqcup\set{0},\, k\bigr)$, then
$p\cOI
p\cong \Mat_k(\cOI )$; \item[(b)] If $\bigl( V(p\cOI
p),\,[p]\bigr)\cong\bigl(\ZZ\sqcup\set{0},\, -k\bigr)$, then
$p\cOI p\cong \Mat_k(p_2\cOI p_2)$; \item[(c)] If $\bigl( V(p\cOI
p),\,[p]\bigr)\cong\bigl(\ZZ\sqcup\set{0},\, 0_{\ZZ }\bigr)$,
then
$p\cOI p\cong p_1\cOI p_1$.
\end{itemize}
\end{remark}

Thus, in order to represent $\ZZ\sqcup\set{0}$ by corners of
$\cOI$, we can restrict our arguments to those corners generated
by $1$, $p_1$, and $p_2$. Furthermore, by
\cite[Theorem~3.5(1)]{LPh}, the isomorphism
\begin{equation}\label{Eq:VMmOn5}
\tau\colon(V(\cOI), [1])\rightarrow (V(p_2\cOI p_2), [p_2])
\end{equation}
is induced by a unital C*-algebra isomorphism
\[
\psi\colon\cOI \rightarrow p_2\cOI p_2.
\]

The remaining basic facts that we shall need are contained in the
following lemmas.

\begin{lemma}\label{L:induceVmapsinf}
Let $B$ be a C*-algebra and let $q\in B$ a projection. Then any
normalized monoid homomorphism
\[
\alpha\colon (V(\cOI),[1]) \rightarrow (V(B),[q])
\]
is induced by a C*-algebra homomorphism
$\phi\colon\cOI\rightarrow
B$ that sends $1$ to $q$. That is, $V(\phi)=\alpha$.
\end{lemma}

\begin{proof}
Set $a=[1]$ and $b_n=[p_n]$, for all $n\geq1$. Observe that
$a=[s_1s_1^*]=a+b_1$ and $b_n=a+b_{n+1}$. As
$[q]=\alpha([1])=[q]+\alpha(b_1)$ and by Lemma~\ref{L:DecKThy},
there exists a projection $q_1\in qBq$ such that $q_1\sim q$ and
$[q-q_1]=\alpha(b_1)$. Suppose that we have constructed pairwise
orthogonal projections $q_1,\dots,q_n\in qBq$ such that $q_l\sim
q$ and $\alpha(b_l)=\bigl[q-\sum_{i=1}^lq_i\bigr]$ for $1\leq
l\leq n$. As
\[
\Bigl[q-\sum_{i=1}^nq_i\Bigr]=\alpha(b_n)=[q]+\alpha(b_{n+1}),
\]
there exists a projection
$q_{n+1}\in\bigl(q-\sum_{i=1}^nq_i\bigr)B\bigl(q-\sum_{i=1}^nq_i\bigr)$
such that $q_{n+1}\sim q$ and
$\bigl[q-\sum_{i=1}^{n+1}q_i\bigr]=\alpha(b_{n+1})$. Thus we have
constructed, by induction, a sequence $\famm{q_i}{i\geq1}$ of
pairwise orthogonal projections in $qBq$ such that $q_n\sim q$
while $\alpha(b_n)=\bigl[q-\sum_{i=1}^nq_i\bigr]$ for all
$n\geq1$. Hence there exists a sequence $\famm{t_i}{i\geq1}$ of
elements of $qBq$ with $t_i^*t_j=q\delta_{i,j}$ and
$t_it_i^*=q_i$
for all positive integers $i,j$. There exists a unique C*-algebra
homomorphism $\phi\colon\cOI\to B$ such that $\phi(1)=q$ while
$\phi(s_i)=t_i$ for all $i\geq1$. As
\begin{align*}
V(\phi)([1])&=[q]=\alpha([1]),\\
V(\phi)([p_n])&=\Bigl[q-\sum_{i=1}^nq_i\Bigr]=\alpha([p_n])
\qquad\text{for all }n\geq1
\end{align*}
and $\set{[1]}\cup\setm{[p_n]}{n\geq1}$ generates the monoid
$V(\cOI)$, we get $V(\phi)=\alpha$.
\end{proof}

\begin{lemma}\label{L:induceVmaps}
Let $A$ be a finite direct sum of full matrix algebras over Cuntz
algebras, let $B$ a C*-algebra, and let $q\in B$ be a projection.
Then any normalized monoid homomorphism
\[
\alpha\colon (V(A),[1_A]) \rightarrow (V(B),[q])
\]
is induced by a C*-algebra homomorphism $\phi\colon A\rightarrow
B$ that sends $1_A$ to $q$. That is, $V(\phi)=\alpha$.
\end{lemma}

\begin{proof}
Write $A= \bigoplus_{j=1}^r\Mat_{k_j}(\cO_{n_j})\oplus
\bigoplus_{i=1}^s \Mat_{l_i}(\cOI)$ for some $k_j,n_j,l_i\in\NN$.
Let~$p$ be the central projection of $A$ corresponding to the
unit
of $\bigoplus_{j=1}^r \Mat_{k_j}(\cO_{n_j})$, and let $q_1$,
\dots
, $q_s$ be the orthogonal central projections in $A$
corresponding
to $\Mat_{l_1}(\cOI)$, \dots, $\Mat_{l_s}(\cOI)$, respectively.
Thus $p$ and $q_1$, \dots , $q_s$ are orthogonal central
projections summing up to $1_A$.

Each $q_i$ is an orthogonal sum of pairwise equivalent
projections
$g_1^{(i)}$, \dots, $g_{l_i}^{(i)}$ such that
$g_1^{(i)}Ag_1^{(i)}
\cong \cOI$. In $V(A)$, we get the equation
\[
[p]+\sum_{i=1}^s l_i[g_1^{(i)}]= [p]+\sum_{i=1}^s\, [q_i]= [1_A],
\]
whence $\alpha([p])+\sum_{i=1}^s l_i\alpha([g_1^{(i)}])= [q]$ in
$V(B)$. By Lemma~\ref{L:DecKThy}, $q$ is an orthogonal sum of
projections $\ol{p},\ol{q}_1,\dots,\ol{q}_s$ of $B$ such that
$\alpha([p])=[\ol{p}]$, $\alpha([q_i])= [\ol{q}_i]$, and
each~$\ol{q}_i$ is an orthogonal sum of pairwise equivalent
projections $h_1^{(i)}$, \dots, $h_{l_i}^{(i)}$ such that
$\alpha([g_1^{(i)}])=[h_1^{(i)}]$.

As $pAp\cong\bigoplus_{j=1}^r \Mat_{k_j}(\cO_{n_j})$, it follows
from \cite[Lemma~7.1]{GPW} that the restriction of $\alpha$ to
$V(pAp)$ is induced by a C*-algebra homomorphism $\phi'\colon
pAp\rightarrow B$ such that $\phi'(p)=\ol{p}$. Furthermore, as
$g_1^{(i)}Ag_1^{(i)}\cong\cOI$, the restriction of~$\alpha$ to
$V(g_1^{(i)}Ag_1^{(i)})$ defines a normalized monoid
morphism~$\alpha_i\colon(V(\cOI), [1])\to(V(B), [h_1^{(i)}])$. By
Lemma~\ref{L:induceVmapsinf}, there exists a C*-algebra
homomorphism $\psi_i\colon g_1^{(i)}Ag_1^{(i)}\rightarrow B$
inducing~$\alpha_i$ such that $\psi_i(g_1^{(i)})=h_1^{(i)}$. As
$q_iAq_i\cong\Mat_{l_i}(g_1^{(i)}Ag_1^{(i)})$ and
$\ol{q}_iB\ol{q}_i\cong\Mat_{l_i}(h_1^{(i)}Bh_1^{(i)})$, the map
$\psi_i$ extends to a C*-algebra homomorphism $\phi_i\colon
q_iAq_i\rightarrow\ol{q}_iB\ol{q}_i$ that induces the restriction
of $\alpha$ to $V(q_iAq_i)$. As the equivalence classes of
projections from $pAp$ and from $g_1^{(i)}Ag_1^{(i)}$, for $1\leq
i\leq s$, generate $V(A)$, the C*-algebra homomorphism
\[
\phi=\phi'\oplus\bigoplus_{i=1}^s\phi_i\colon
A\to\ol{p}B\ol{p}\oplus\bigoplus_{i=1}^s\ol{q}_iB\ol{q}_i=qBq\subseteq
B
\]
induces $\alpha$.
\end{proof}

\begin{theorem}\label{T:Cstarnonun} An abelian monoid $M$
is isomorphic to $V(A)$ for some special \pup{resp., extended}
Cuntz limit~$A$ if and only if
\begin{itemize}
\item[(a)] $M$ is a countable, regular conical refinement monoid.
\item[(b)] For all idempotents $e\leq f$ in $M$, the homomorphism
$\AC{M}{e} \rightarrow \AC{M}{f}$ given by $x\mapsto x+f$ is
injective, and $\AC{M}{e}+f$ is a pure subgroup of $\AC{M}{f}$.
\end{itemize}
\end{theorem}

\begin{proof}
$(\Longrightarrow)$: By \eqref{Eq:VMmOn} and \eqref{Eq:VMmOn2}
and
since the $V(-)$ functor preserves direct limits and finite
direct
sums, the present implication follows from the easy direction of
Theorem~\ref{T:oL=oR}.

$(\Longleftarrow)$: Since $M$ is countable, Theorem~\ref{T:oL=oR}
implies that $M$ is the direct limit of a sequence of the form
\[
M_1\xrightarrow{\alpha_1}M_2\xrightarrow{\alpha_2}
M_3\xrightarrow{\alpha_3}\cdots
\]
where each $M_i$ is a finite direct sum of monoids
$\ZZ\sqcup\set{0}$ and $(\ZZ/n_{ij}\ZZ)\sqcup\set{0}$ for some
$2\leq n_{ij}<\infty$. Hence, denoting by $A_i$ the direct sum of
the corresponding Cuntz algebras $\cOI$ and $\cO_{n_{ij}+1}$,
there is an isomorphism $h_i\colon V(A_i)\rightarrow M_i$. Each
of
the homomorphisms
\[
h_{i+1}^{-1}\alpha_ih_i\colon V(A_i) \longrightarrow V(A_{i+1})
\]
sends $[1_{A_i}]$ to the class of a projection in $A_{i+1}$, and
so, by Lemma~\ref{L:induceVmaps}, $h_{i+1}^{-1}\alpha_ih_i$ is
induced by a C*-algebra homomorphism $\phi_i\colon A_i\rightarrow
A_{i+1}$. Therefore $M \cong V(A)$, where $A$ is the C*-inductive
limit of the sequence
\begin{equation*}
A_1\xrightarrow{\phi_1}A_2\xrightarrow{\phi_2}
A_3\xrightarrow{\phi_3}\cdots\tag*{\qed}
\end{equation*}
\renewcommand{\qed}{}
\end{proof}

A structural description of the  monoids appearing in
Theorem~\ref{T:Cstarnonun} is easily obtained with the help of
\cite[Theorem 3.3]{GPW}, as follows.

\begin{corollary}\label{C:VAnonun}
Let $M$ be an abelian monoid. Then $M\cong V(A)$ for some special
\pup{resp., extended} Cuntz limit $A$ if and only if
\[
M \cong \bigsqcup_{e\in\Lambda} \bigl( \set{e}\times G_e \bigr)
\subseteq \Lambda\times G
\]
where
\begin{itemize}
\item[(a)] $\Lambda$ is a countable distributive semilattice.
\item[(b)] $G$ is a countable abelian group. \item[(c)] $G_e$ is
a
pure subgroup of $G$ for all $e\in\Lambda$. \item[(d)]
$G_0=\set{0}$ and $\bigcup_{e\in\Lambda} G_e =G$. \item[(e)]
$G_e+G_f= G_{e+f}$ and $G_e\cap G_f= \bigcup_{g\in\Lambda,\,
g\leq
e,f} G_g$ for all $e$, $f\in\Lambda$.
\end{itemize}
\end{corollary}

\begin{corollary}\label{C:FromExt2Sp}
For every extended Cuntz limit~$A$, there exists a special Cuntz
limit~$B$ such that $V(A)\cong V(B)$.
\end{corollary}

\section{Algebras with order-unit}\label{S:OrdUnits}

In the present section, we establish unital versions of
Theorem~\ref{T:Cstarnonun} and Corollary~\ref{C:VAnonun}. For
this, we need suitable adaptations of
Lemma~\ref{L:induceVmapsinf}
to the corners of~$\cOI$ generated by the projections $p_1$ and
$p_2$.

\begin{lemma}\label{L:induceVmapsinfP2}
Let $B$ be a C*-algebra, and $q\in B$ a projection. Then any
normalized monoid homomorphism
\[
\alpha\colon (V(p_2\cOI p_2),[p_2]) \rightarrow (V(B),[q])
\]
is induced by a C*-algebra homomorphism $\phi\colon p_2\cOI
p_2\rightarrow B$ that sends $p_2$ to $q$. That is,
$V(\phi)=\alpha$.
\end{lemma}

\begin{proof}
As noticed in \eqref{Eq:VMmOn5}, the normalized monoid
isomorphism
\[
\tau\colon(V(\cOI), [1])\rightarrow (V(p_2\cOI p_2))
\]
is induced by a unital C*-algebra isomorphism $\psi\colon\cOI
\rightarrow p_2\cOI p_2$, that is, $\tau =V(\psi)$. Thus, $\beta
=\alpha \tau\colon (V(\cOI),[1]) \rightarrow (V(B),[q])$ is a
normalized monoid morphism. By Lemma~\ref{L:induceVmapsinf},
there
exists a C*-algebra homomorphism $\varphi\colon\cOI\rightarrow B$
sending $1$ to $q$, such that $V(\varphi )=\beta =\alpha \tau$.
Since $\psi $ is an isomorphism and $V(-)$ is a functor, we get
$\tau ^{-1}=V(\psi ^{-1})$, and thus $\alpha =V(\varphi \psi
^{-1})$, where $\varphi\psi^{-1}\colon p_2\cOI p_2\rightarrow B$
is a C*-algebra homomorphism that sends $p_2$ to $q$. Thus,
$\phi=\varphi\psi^{-1}$ is the desired morphism.
\end{proof}

In the case of the corner $p_1\cO_{\!\infty\,}p_1$, as $[p_1]$ is
idempotent, we need to restrict the target algebras in order to
preserve the ``lifting'' property.

\begin{lemma}\label{L:induceVmapsinfP1}
Let $B$ be a finite direct sum of full matrix algebras over Cuntz
algebras~$\cO_n$ \pup{for $2\leq n<\infty$} and $p\cOI p$
\pup{for
any nonzero projection $p\in \cOI$}, and let $q\in B$ be a
projection. Then any normalized monoid homomorphism
\[
\alpha\colon (V(p_1\cOI p_1),[p_1]) \rightarrow (V(B),[q])
\]
is induced by a C*-algebra homomorphism $\phi\colon p_1\cOI
p_1\rightarrow B$ that sends $p_1$ to $q$. That is,
$V(\phi)=\alpha$.
\end{lemma}

\begin{proof}
By Remark~\ref{Rk:Classification2}, we can write
$B=\bigoplus_{l=1}^{r+s}B_l$, where $B_j=\Mat_{k_j}(\cO_{n_j})$
for $1\leq j\leq r$ (with $k_j,n_j\in\NN$) and
$B_{r+i}=p_{l_i}\cOI p_{l_i}$ for $1\leq i\leq s$ (with
$l_i\in\ZZ^+$). There exist pairwise orthogonal projections
$q_i\in B_i$ such that $q=\sum\limits_{i=1}^{r+s}q_i$. Since the
functor $V(-)$ preserves finite direct sums, we can reduce the
problem to the case where~$B$ is either $\Mat_{k_j}(\cO_{n_j})$
or
$p_{l_i}\cOI p_{l_i}$ with $l_i\in\ZZ^+$, by composing with the
canonical projections $\pi _i\colon B\onto B_i$. We get
$\alpha_i=V(\pi_i)\alpha\colon V(p_1\cOI p_1) \rightarrow V(B_i)$
with $\alpha _i([p_1])=[q_i]$. If $q_i=0$, then $\alpha_i$ is
induced by the zero homomorphism $\phi_i\colon p_1\cOI p_1\to
B_i$; so suppose that $q_i\neq0$. Since $[p_1]$ is an order-unit
of $V(p_1\cOI p_1)$, the image of $\alpha _i$ is contained in the
order-ideal of $V(B_i)$ generated by $[q_i]$, whence~$\alpha _i$
restricts to a normalized monoid morphism
\[
\alpha_i\colon (V(p_1\cOI p_1),[p_1]) \rightarrow
(V(q_iB_iq_i),[q_i])
\]
Since $B_i$ is a purely infinite simple C*-algebra, so is
$q_iB_iq_i$.

Since $V(p_1\cOI p_1)\setminus\set{0}$ is a group containing
$[p_1]$, $[q_i]\ne 0$ in $V(q_iB_iq_i)$, and $V(q_iB_iq_i)$ is a
conical monoid, it follows from \cite[Corollary~2.2]{AGP} that
\[
\beta_i={\alpha _i}\res_{V(p_1\cOI p_1)\setminus\set{0}}\colon
K_0(p_1\cOI p_1)\rightarrow K_0(q_iB_iq_i)
\]
is a group homomorphism such that $\beta_i([p_1])=[q_i]$. By
Remark~\ref{Rk:Classification}, for each $i\geq 1$, $q_iB_iq_i$
is
isomorphic to either $\Mat_k(\cO_n )$ (for some $k\geq 1$) or
$p\cOI p$ (for some projection $p\in \cOI$). Thus, by
\cite[Lemma~3.7]{LPh}, there exists a unital C*-algebra
homomorphism $\phi_i\colon p_1\cOI p_1\rightarrow q_iB_iq_i$ such
that $K_0(\phi_i)=\beta_i$, and thus $V(\phi_i)=\alpha_i$. The
map
\[
\phi=\bigoplus_{i=1}^{r+s}\phi_i\colon p_1\cOI p_1\rightarrow
qBq\subseteq B
\]
satisfies the desired properties.
\end{proof}

Thus, we get the following version of Lemma~\ref{L:induceVmaps}.

\begin{lemma}\label{L:induceVmapsUnit}
Let $A,B$ be finite direct sums of full matrix algebras over
either Cuntz algebras $\cO_n$ \pup{$2\leq n<\infty$} or $p\cOI p$
\pup{for projections $p\in\cOI$}. Then any normalized monoid
homomorphism
\[
\alpha\colon (V(A),[1_A]) \rightarrow (V(B),[1_B])
\]
is induced by a C*-algebra homomorphism $\phi\colon A\rightarrow
B$ that sends $1_A$ to $1_B$. That is, $V(\phi)=\alpha$.
\end{lemma}

\begin{proof}[Outline of proof]
By arguing as in the proof of Lemma~\ref{L:induceVmaps} and using
Remark~\ref{Rk:Classification2} together with \eqref{Eq:VMmOn3},
we reduce the problem to the case where $A$ is either $\cOI$, or
$p_1\cOI p_1$, or $p_2\cOI p_2$. The first case is covered by
Lemma~\ref{L:induceVmapsinf}, the second case by
Lemma~\ref{L:induceVmapsinfP1}, and the third case by
Lemma~\ref{L:induceVmapsinfP2}.
\end{proof}

\begin{theorem}\label{T:seastar} Let $(M,u)$ be an abelian monoid
with order-unit. Then $(M,u)\cong (V(A),[1_A])$ for some unital
extended Cuntz limit $A$ if and only if
\begin{itemize}
\item[(a)] $M$ is a countable, regular conical refinement monoid.
\item[(b)] For all idempotents $e\leq f$ in $M$, the homomorphism
$\AC{M}{e} \rightarrow \AC{M}{f}$ given by $x\mapsto x+f$ is
injective, and $\AC{M}{e}+f$ is a pure subgroup of $\AC{M}{f}$.
\end{itemize}
\end{theorem}
\begin{proof} $(\Longrightarrow)$: Theorem~\ref{T:Cstarnonun}.

$(\Longleftarrow)$: Corollary~\ref{C:Limorderunit} implies that
$(M,u)$ is the direct limit of a sequence of the form
\[
(M_1,u_1)\xrightarrow{\alpha_1}(M_2,u_2)\xrightarrow{\alpha_2}
(M_3,u_3)\xrightarrow{\alpha_3}\cdots
\]
where each $(M_i,u_i)$ is a finite direct sum of pairs
$((\ZZ/n_{ij}\ZZ)\sqcup\set{0},\, \mbar_{ij})$ and
$(\ZZ\sqcup\set{0}, k_i)$ for some $n_{ij},m_{ij}\in\NN$, $k_i\in
\ZZ$. In view of~\eqref{Eq:VMmOn}--\eqref{Eq:VMmOn3}, there exist
isomorphisms $h_i\colon (V(A_i),[1_{A_i}]) \rightarrow (M_i,u_i)$
where $A_i$ is a direct sum of matrix algebras of the form either
$\Mat_{m_{ij}}(\mathcal{O}_{n_{ij}+1})$ or $p_{k_i}\cOI p_{k_i}$,
with $p_{k_i}$ being suitable projections. By
Lemma~\ref{L:induceVmapsUnit}, each of the normalized
homomorphisms
\[
h_{i+1}^{-1}\alpha_ih_i\colon (V(A_i),[1_{A_i}]) \longrightarrow
(V(A_{i+1}),[1_{A_{i+1}}])
\]
is induced by a unital C*-algebra homomorphism $\phi_i\colon
A_i\rightarrow A_{i+1}$. Therefore $(M,u)\cong(V(A),[1_A])$ where
$A$ is the C*-inductive limit of the sequence
\begin{equation*}
A_1\xrightarrow{\phi_1}A_2\xrightarrow{\phi_2}A_3\xrightarrow{\phi_3}
\cdots\tag*{\qed}
\end{equation*}
\renewcommand{\qed}{}
\end{proof}

\begin{corollary}\label{C:VAstruct}
Let $(M,u)$ be an abelian monoid with order-unit. Then
$(M,u)\cong
(V(A),[1_A])$ for some unital extended Cuntz limit $A$ if and
only
if
\[
(M,u) \cong \Bigl( \bigsqcup_{e\in\Lambda} \bigl( \set{e}\times
G_e \bigr),\, (1,u_1) \Bigr) \subseteq \bigl( \Lambda\times
G_1,\,
(1,u_1) \bigr)
\]
where
\begin{itemize}
\item[(a)] $\Lambda$ is a countable distributive semilattice with
maximum element $1$. \item[(b)] $G_1$ is a countable abelian
group. \item[(c)] $G_e$ is a pure subgroup of $G_1$ for all
$e\in\Lambda$, and $G_0=\set{0}$. \item[(d)] $G_e+G_f= G_{e+f}$
and $G_e\cap G_f= \bigcup_{g\in\Lambda,\, g\leq e,f} G_g$ for all
$e$, $f\in\Lambda$. \item[(e)] $u_1\in G_1$.
\end{itemize}
\end{corollary}
\begin{proof}
$(\Longrightarrow$): By Corollary~\ref{C:VAnonun}, $M$ is
isomorphic to a monoid of the form
\[
M'= \bigsqcup_{e\in\Lambda} \bigl( \set{e}\times G_e \bigr)
\subseteq \Lambda\times G
\]
for some countable distributive semilattice $\Lambda$ and some
countable abelian group~$G$ with subgroups $G_e$ satisfying the
conditions of that corollary. As $M$ has an order-unit, $\Lambda$
has a largest element. Conditions (a)--(e) are now all satisfied.

$(\Longleftarrow)$: With the help of \cite[Theorem 3.3]{GPW}, it
is clear that $M$ satisfies conditions (a) and (b) of
Theorem~\ref{T:seastar}.
\end{proof}

\section*{Acknowledgments}

Part of this work was done during a visit of the second author to
the Departamento de Matem\'aticas de la Universidad de C\'adiz.
The second author wants to thank the host center for its warm
hospitality. Both authors thank Pere Ara for turning our attention
to reference \cite{AF}, and also thank the referees for their
comments, that really improved the final version of this paper.

\end{document}